\documentclass{elsart3-1}

\usepackage{amsfonts}
\usepackage{amsmath}
\usepackage{amssymb}
\usepackage[english,francais]{babel}
\def\Z{\mathbb{Z}}
\def\R{\mathbb{R}}

\def\epsilon{\varepsilon}
\newtheorem{theorem}{Theorem}[section]
\newtheorem{lemma}[theorem]{Lemma}
\newtheorem{e-proposition}[theorem]{Proposition}
\newtheorem{corollary}[theorem]{Corollary}
\newtheorem{e-definition}[theorem]{Definition\rm}
\newtheorem{remark}{\bf Remark\/}

\newtheorem{theoreme}{Th\'eor\`eme}[section]
\newtheorem{lemme}[theoreme]{Lemme}
\newtheorem{proposition}[theoreme]{Proposition}

\renewcommand{\theequation}{\arabic{equation}}
\def\ds{\displaystyle}
\def\og{\leavevmode\raise.3ex\hbox{$\scriptscriptstyle\langle\!\langle$~}}
\def\fg{\leavevmode\raise.3ex\hbox{~$\!\scriptscriptstyle\,\rangle\!\rangle$}}

\begin{document}

\begin{frontmatter}


\selectlanguage{english}
\title{Asymptotics of the  KPP  minimal speed within large drift}


\selectlanguage{francais}
\author[authorlabel1]{Mohammad El Smaily} \ead{elsmaily@math.ubc.ca} and
\author[authorlabel1]{St\'ephane kirsch}
\ead{kirsch@math.ubc.ca}

\address[authorlabel1]{{\footnotesize Department of Mathematics, University of British Columbia }{\footnotesize  $\&$ Pacific Institute for the Mathematical Sciences}
 \\{\footnotesize 1984 Mathematics Road, V6T 1Z2, Vancouver, BC, Canada}}

\selectlanguage{english}
\begin{abstract}
This Note is concerned with the asymptotic behavior of the minimal KPP speed of propagation for reaction-advection-diffusion equations with a large
drift $Mq$ (where $q$ is the advection).  We first give the limit of the speed as
$M\rightarrow+\infty$ in any space dimension $N.$ Then, we give the necessary and sufficient condition that the advection field should satisfy  so that the speed  acts as $O(M)$ as $M\rightarrow+\infty.$
\vskip 0.5\baselineskip
\selectlanguage{francais}
\noindent{\bf R\'esum\'e}
\vskip 0.5\baselineskip
\noindent Cette Note est consacr\'ee \`{a} l'\'{e}tude asymptotique de la vitesse minimale de propagation pour les fronts progressifs pulsatoires qui
satisfont une \'{e}quation de r\'{e}action-advection-diffusion lors d'une grande advection $Mq$ (o\`{u} $q$ est l'advection).  Nous donnons la limite de la vitesse quand $M\rightarrow+\infty$ dans une dimension $N$ quelconque. Puis, quand $N=2,$ on donne une condition n\'{e}cessaire et suffisante pour que la vitesse  se comporte comme $O(M)$
lorsque $M\rightarrow+\infty.$
\end{abstract}

\end{frontmatter}
\selectlanguage{francais}
\section*{Version fran\c{c}aise abr\'eg\'ee}
Pour chaque $M>0,$ nous consid\'{e}rons l'\'{e}quation de r\'{e}action-advection-diffusion suivante
\begin{equation}\label{het eq}
 \left\{
      \begin{array}{l}
 u_t =\nabla\cdot(A(z)\nabla u)\;+M\,q(z)\cdot\nabla u+f(z,u),\; t\in\,\mathbb{R},\;z\in\,\Omega,
\vspace{3pt} \\
        \nu \cdot A\nabla u =0\;\hbox{ sur }\mathbb{R}\times\partial\Omega,
      \end{array}
    \right.
\end{equation}
o\`{u} $\nu$ est la normale unitaire ext\'erieure sur $\partial\Omega$ quand $\partial\Omega\neq\emptyset.$ Le domaine $\Omega$ est un sous-ensemble connexe de $\mathbb{R}^{N}$ de classe $C^3$ pour lequel il existe $1\leq d\leq N,$ $L_1,\cdots,L_d$ positifs et $R>0$ tels que
$$\forall\,(x,y)\,\in\,\Omega\subseteq\R^{d}\times\R^{N-d},\,|y|\,\leq\,R,\hbox{ et }\forall\,k=(k_1,\cdots,k_d,0,\cdots,0)\in L_1\mathbb{Z}\times\cdots\,\times L_d\mathbb{Z}\times\{0\}^{N-d},\displaystyle{\Omega=\Omega+k}.$$ On note la cellule de p\'{e}riodicit\'e de $\Omega$ par
        $$C=\{(x,y)\in\,\Omega
;\; x_{1}\in(0,L_1),\cdots,x_{d}\in(0,L_d)\}.$$ Dans ce cadre p\'{e}riodique, un champ $v:\Omega\rightarrow\,\mathbb{R}^N$ est dit $L$-periodique en  $x$
 si $w(x+k,y)=w(x,y)$ p.p. dans
$ \Omega$ quel que soit $\displaystyle{k=(k_1,\cdots\,,k_d)\in\prod^{d}_{i=1}L_i\mathbb{Z}}.$ La diffusion $A(x,y)=(A_{ij}(x,y))_{1\leq i,j\leq
 N}$ dans l'\'{e}quation (\ref{het eq}) est un champ matriciel de classe $C^{2,\delta}(\,\overline{\Omega}\,)$ (avec $\delta
\,>\,0$) qui est $L$-p\'{e}riodique en $x$ et v\'erifie
        $$\exists\,0<\alpha_1\leq\alpha_2,\forall(x,y)\;\in\;\Omega,\forall\,\xi\,\in\,\mathbb{R}^N,
       \displaystyle{ \alpha_1|\xi|^2 \;\leq\;\sum_{1\leq i,j\leq N}\,A_{ij}(x,y)\xi_i\xi_j\,\;\leq\alpha_2|\xi|^2.}$$
 L'advection $q(x,y)=(q_1(x,y),\cdots,q_N(x,y))$  est un champ vectoriel $L$-p\'{e}riodique en $x$ et de classe
$C^{1,\delta}(\overline{\Omega})$ ($\delta>0$) qui satisfait
        \begin{equation}\label{cqf}
        \nabla\cdot q=0\,\hbox{ dans }\, \overline{\Omega},\;
         q\cdot\nu=0\; \hbox{sur}\;\partial\Omega\hbox{ (quand $\partial\Omega\neq\emptyset$)},\hbox{ et }
        \forall\,1\leq i\leq d, \;\displaystyle{\int_{C}q_i\;dx\,dy =0}\hbox{.}
        \end{equation}
La partie non lin\'{e}are  $f=f(x,y,u)$ est une fonction positive de classe $C^{1,\delta}(\overline{\Omega}\times[0,1]),$ $L$-p\'{e}riodique en $x$  telle que $ \displaystyle{f(x,y,0)=f(x,y,1)=0 }$ pour tout
$(x,y)\in\overline{\Omega}$ et
\begin{eqnarray}\label{cf1}
    \left\{
      \begin{array}{ll}
        \exists \,\rho\in(0,1),\;\forall(x,y)\,\in\overline{\Omega},\;\displaystyle{\forall\, 1-\rho\leq s \leq
s'\leq1,}\;
\displaystyle{f(x,y,s)\;\geq\,f(x,y,s')} \hbox{,} \vspace{3 pt}\\
        \forall\,(x,y)\in\overline{\Omega},\quad \zeta(x,y):=\displaystyle{f'_{u}(x,y,0)=\lim_{u\rightarrow\,0^+}\frac{f(x,y,u)}{u}>0}  \hbox{.}
      \end{array}
    \right.
\end{eqnarray}
 On suppose aussi que la r\'eaction $f$ satisfait la condition  ``KPP'' (d'apr\`{e}s Kolmogorov, Petrovsky et Piskunov \cite{KPP})
 \begin{equation}\label{cf2}
 \forall\, (x,y,s)\in\overline{\Omega}\times(0,1),~0<f(x,y,s)\leq f'_u(x,y,0)\,s=\zeta(x,y)s.
\end{equation}
Un exemple de cette non lin\'{e}arit\'{e} est la fonction ``homog\`{e}ne'' $f(u)=u(1-u)$ sur $(0,1).$

On s'int\'eresse au ph\'{e}nom\`{e}ne de propagation des fronts progressifs pulsatoires pour l'\'{e}quation (\ref{het eq}). Nous fixons une
 direction unitaire $e\in\mathbb{R}^{d}$ et nous notons  $\tilde{e}=(e,0,\cdots,0)\in\mathbb{R}^N.$ Un front progressif
pulsatoire qui se propage \emph{dans la direction de $-e$ avec une vitesse }$c$ est une solution $u(t,x,y)$ de (\ref{het eq}) de la forme $u(t,x,y)=\phi(x\cdot e+ct,x,y)$ o\`{u} la fonction $\phi$ est $L$-p\'{e}riodique en $x$ et satisfait les conditions limites $\ds{\lim_{s\rightarrow-\infty}\phi(s,x,y)=0\hbox{ et }\lim_{s\rightarrow+\infty}\phi(s,x,y)=1}$ uniform\'ement en $\ds{(x,y)\in\overline{\Omega}}.$  \emph{D'apr\`{e}s les r\'{e}sultats de \cite{bh} et \cite{bhn1}, il existe une valeur critique not\'{e}e $\ds{c^{*}_{\Omega,A,Mq,f}(e)}>0,$ qui s'appelle  la vitesse minimale KPP, telle qu'il existe un front progressif pulsatoire pour l'\'{e}quation (\ref{het eq}) avec une vitesse $c$ si et seulement si $c\geq \ds{c^{*}_{\Omega,A,Mq,f}(e)}>0.$} De plus, selon \cite{bhn2}, le terme
 param\'{e}trique $\ds{c^{*}_{\Omega,A,Mq,f}(e)}/M$ reste born\'{e} ind\'{e}pendamment de $M\geq1.$ Nous avons d\'{e}termin\'e, dans \cite{ek1}, la limite de $\ds{c^{*}_{\Omega,A,Mq,f}(e)}/M$ lorsque $M\rightarrow+\infty.$
  Dans le th\'{e}or\`{e}me suivant, nous donnons des d\'{e}tails sur la  limite $\ds{c^{*}_{\Omega,A,Mq,f}(e)}/M$ lors que la dimension est $N=2.$

\begin{lemme} [trajectoires p\'{e}riodique non born\'{e}es de $q$] On suppose ici que la dimension est $N=2$ et donc $d\in\{1,2\}.$
Soit $T(x)$ une  trajectoire p\'{e}riodique non born\'{e}e de $q$ dans $\Omega$ passant par $x.$ C'est \`{a} dire, il existe $\mathbf{a}\in L_1\Z\times L_2\Z\setminus\{0\}$ (resp. $L_1\Z\times \{0\}\setminus\{0\}$) quand $d=2$ (resp. $d=1$) tel que $T(x)=T(x)+\mathbf{a}.$ Dans ce cas, on dit $T(x)$ est $\mathbf{a}-$p\'{e}riodique. Alors, toutes les autres trajectoires non born\'{e}es p\'{e}riodiques $T(y)$ de $q$ sont aussi $\mathbf{a}-$p\'{e}riodique. De plus, $\mathbf{a}=L_1e_1$ quand $d=1$ et donc, dans le cas $d=1,$  toutes les trajectoires p\'{e}riodiques sont
$L_1e_1-$ p\'{e}riodiques.
\end{lemme}

\vskip0.2cm

\begin{theoreme} Supposons que $N=2$ et que $q,$ $\Omega,$ $A$ et $f$ satisfont les hypoth\`{e}ses mentionn\'{e}s au dessus (avec $N=2$). Alors, \\
(i) S'il n'existe pas de trajectoire p\'{e}riodique non born\'{e}e de $q,$ alors
$\displaystyle{
\lim_{M\rightarrow+\infty}\ds{\frac{\ds{c^{*}_{\Omega,A,M\,q,f}(e)}}{M}}=0},
$ quelle que soit la direction unitaire $e.$\\
(ii) S'il existe une trajectoire p\'{e}riodique non born\'{e}e  $T(x)$ de $q$ dans $\Omega$ (\`a la quelle on peut donc associer une p\'eriode $\mathbf{a}\in\mathbb{R}^2$), alors \begin{equation}\label{null equivalence}
\lim_{M\rightarrow+\infty}\ds{\frac{\ds{c^{*}_{\Omega,A,M\,q,f}(e)}}{M}} >0\,\Longleftrightarrow\,\tilde{e}\cdot\mathbf{a}\neq0.
\end{equation}
De plus, dans le cas o\`{u} $d=1,$ le Lemme \ref{periodic trajs} implique que  $\tilde{e}\cdot\mathbf{a}=\pm L_1\neq0.$ Utilisant (\ref{null equivalence}), on peut alors conclure que, pour $d=1,$
\begin{equation}\label{null equivalence d=1}
\lim_{M\rightarrow+\infty}\ds{\frac{\ds{c^{*}_{\Omega,A,M\,q,f}(e)}}{M}} >0\,\Longleftrightarrow\,\text{( Il existe une trajectoire p\'{e}riodique non born\'{e}e $T(x)$ de $q$ dans $\Omega$)}.
\end{equation}

\end{theoreme}

\selectlanguage{english}
\section{Asymptotics of the minimal speed within large drift in any dimension $N$ with more details in the case $N=2$}
For each $M>0,$  we consider the reaction-advection-diffusion equation
\begin{equation}\label{het eq}
 \left\{
      \begin{array}{l}
 u_t =\nabla\cdot(A(z)\nabla u)\;+M\,q(z)\cdot\nabla u+f(z,u),\; t\in\,\mathbb{R},\;z\in\,\Omega,
\vspace{3pt} \\
        \nu \cdot A\nabla u =0\;\hbox{ on }\mathbb{R}\times\partial\Omega,
      \end{array}
    \right.
\end{equation}
where $\nu$ stands
for the unit outward normal on $\partial\Omega$ whenever it is nonempty.

The domain $\Omega$ is a
$C^3$ nonempty connected open subset of $\mathbb{R}^N$ such that for some integer $1\leq d\leq N,$ there exist $L_1,\cdots,L_d$  positive real numbers such that
\begin{eqnarray}\label{comega}
    \left\{
      \begin{array}{l}
        \exists\,R\geq0\,;\forall\,(x,y)\,\in\,\Omega\subseteq\R^{d}\times\R^{N-d},\,|y|\,\leq\,R, \\
        \forall\,(k_1,\cdots,k_d)\in\,L_1\mathbb{Z}\times\cdots\,\times L_d\mathbb{Z},
        \quad\displaystyle{\Omega\;=\;\Omega+\sum^{d}_{k=1}k_ie_i},
      \end{array}
    \right.
\end{eqnarray}
where $\;(e_i)_{1\leq i\leq N}\;$ is the canonical basis of
$\mathbb{R}^N.$ In other words, $\Omega$ is bounded in the $y-$direction and periodic in $x.$ As archetypes of the domain $\Omega,$ we may have the whole space $\R^{N}$ which corresponds for $d=N$ and $L_1,\cdots, L_N$ any array of positive real numbers. We may also have the whole space $\R^N$ with a periodic
array of holes or an infinite cylinder with an oscillating boundary.
In this periodic situation, we call
\begin{equation}\label{periodicity cell}
C=\{(x,y)\in\,\Omega
;\; x_{1}\in(0,L_1),\cdots,x_{d}\in(0,L_d)\}
\end{equation}  the
periodicity cell of $\Omega.$

The diffusion matrix $A(x,y)=(A_{ij}(x,y))_{1\leq i,j\leq
 N}$ is a \textit{symmetric} $C^{2,\delta}(\,\overline{\Omega}\,)$ (with $\delta
\,>\,0$) matrix field which is $L$-periodic with respect to $x$ and satisfies
\begin{eqnarray}\label{cA}
      \begin{array}{l}
         \exists\,0<\alpha_1\leq\alpha_2,\forall(x,y)\;\in\;\Omega,\forall\,\xi\,\in\,\mathbb{R}^N,\hbox{ we have }
      \; \displaystyle{ \alpha_1|\xi|^2 \;\leq\sum_{1\leq i,j\leq N}A_{ij}(x,y)\xi_i\xi_j\,\;\leq\alpha_2|\xi|^2.}
      \end{array}
\end{eqnarray}

The underlying advection $q(x,y)=(q_1(x,y),\cdots,q_N(x,y))$  is a
$C^{1,\delta}(\overline{\Omega})$ (with $\delta>0$) vector field which is  $L$- periodic with respect to $x$
and satisfies
\begin{equation}\label{cq}
      \begin{array}{llll}
        \nabla\cdot q=0\hbox{ in } \overline{\Omega}, &
         q\cdot\nu=0\hbox{ on }\partial\Omega\hbox{ (when $\partial\Omega\neq\emptyset$)}, &
        \forall\,1\leq i\leq d,\displaystyle{\int_{C}q_i\;dx\,dy =0}\hbox{.}
      \end{array}
\end{equation}

Concerning the nonlinearity  $f=f(x,y,u),$ it is a nonnegative
function defined in $\overline{\Omega}\,\times[0,1],\;$ such that
\begin{eqnarray}\label{cf1}
    \left\{
      \begin{array}{ll}
       f\;\hbox{ is $L$-periodic with respect to }\; x, \hbox{ and of class }C^{1,\delta}(\overline{\Omega}\times[0,1]),\vspace{3 pt}\\
        \forall\,(x,y)\in\,\overline{\Omega},\quad \displaystyle{f(x,y,0)=f(x,y,1)=0 } \hbox{,} \vspace{3 pt}\\
        \exists \,\rho\in(0,1),\;\forall(x,y)\,\in\overline{\Omega},\;\displaystyle{\forall\, 1-\rho\leq s \leq
s'\leq1,}\;
\displaystyle{f(x,y,s)\;\geq\,f(x,y,s')} \hbox{,} \vspace{3 pt}\\
         \forall\,(x,y)\in\overline{\Omega},\quad \zeta(x,y):=\displaystyle{f'_{u}(x,y,0)=\lim_{u\rightarrow\,0^+}\frac{f(x,y,u)}{u}>0}  \hbox{,}
      \end{array}
    \right.
\end{eqnarray}
 with the additional ``KPP'' assumption (referring  to \cite{KPP} by Kolmogorov, Petrovsky and  Piskunov)
 \begin{equation}\label{cf2}
 \forall\, (x,y,s)\in\overline{\Omega}\times(0,1),~0<f(x,y,s)\leq f'_u(x,y,0)\times\,s.
\end{equation}
An archetype of $f$ is $(x,y,u)\mapsto u(1-u)h(x,y)$ defined on $\overline{\Omega}\times[0,1]$ where $h$ is a positive $C^{1,\delta}(\,\overline{\Omega}\,)$ $L$-periodic function.

In all of this paper, $e\in\R^{d}$ is a fixed unit vector and $\tilde{e}:=(e,0,\cdots,0)\in\R^{N}.$ A pulsating travelling front propagating in the direction of $-e$ within a speed $c\neq0$ is a solution $u=u(t,x,y)$ of (\ref{het eq}) for which there exists a function $\phi$ such that
$u(t,x,y)=\phi(x\cdot e+ct,x,y),$  $\phi$ is $L$-periodic in $x$ and $$\lim_{s\rightarrow-\infty}\phi(s,x,y)=0\hbox{ and }\lim_{s\rightarrow+\infty}\phi(s,x,y)=1,$$ uniformly in $(x,y)\in\overline{\Omega}.$

In the same setting as in this paper, it was proved in \cite{bh} and \cite{bhn2} that for all $\Omega,$ $A,\,q,$ and $f$ satisfying (\ref{comega}), (\ref{cA}), (\ref{cq}), and (\ref{cf1}) respectively, there exists $\ds{c^{*}_{\Omega,A,q,f}(e)},$ called the \emph{minimal speed of propagation}, such that pulsating travelling fronts exist if and only if $c\geq \ds{c^{*}_{\Omega,A,q,f}(e)}.$ This result extended that of \cite{KPP} which proved that $c^*(e)=2\sqrt{f'(0)}$ in a ``homogenous'' framework where $f=f(u)$ and there is no advection $q.$ A variational formula for the minimal speed $\ds{c^{*}_{\Omega,A,q,f}(e)}$ involving the principal eigenvalue of an elliptic operator was proved in \cite{bhn2}. Moreover,  El Smaily \cite{El Smaily min max} proved a $\min$-$\max$ formula for the minimal speed. In the following, we recall the definition of ``first integrals'' of
a vector field which was introduced in \cite{bhn1}.

\vskip 0.3cm

\begin{e-definition}[First integrals]\label{first integral} The family of first integrals of an incompressible advection $q$ of the type (\ref{cq}) is defined by
\begin{equation*}
\begin{array}{ll}
\mathcal{I}:=&\left\{w\in H^{1}_{loc}(\overline{\Omega}),\,w\neq0,\;w \hbox{ is } L-\hbox{periodic in $x,$ and }q\cdot\nabla w=0\hbox{ almost everywhere in }\Omega\right\}.
\end{array}
\end{equation*}
Having a matrix $A=A(x,y)$ of the type (\ref{cA}), we also define
\begin{equation}\label{I1}
\mathcal{I}_1^A:=\left\{w\in\,\mathcal{I},\hbox{ such that }\int_C\zeta w^2\geq\int_C\nabla w\cdot A\nabla w\right\}.
\end{equation}
\end{e-definition}

The following theorem gives the asymptotic behavior of the minimal speed in the presence of a large advection in any dimension $N$.
\vskip0.3cm
\begin{theorem}\label{main}
We fix a unit direction $e\in\R^d$ and assume that the diffusion matrix $A$ and the nonlinearity $f$ satisfy (\ref{cA}), (\ref{cf1}) and (\ref{cf2}).
Let $q$ be an advection field which satisfies (\ref{cq}). Then,
\begin{equation}\label{large advection}
\lim_{M\rightarrow+\infty}\ds{\frac{\ds{c^{*}_{\Omega,A,M\,q,f}(e)}}{M}}=\ds\max_{\ds{w\in \mathcal{I}_1^A}}\frac{\ds\int_C(q\cdot\tilde{e})\,w^2}{\ds\int_C w^2}.
\end{equation}
\end{theorem}
The above theorem was proved in details in \cite{ek1} and \cite{zlatos}. In \cite{zlatos}, Zlato\v{s} treated the problem when the domain $\Omega$ is the whole $\mathbb{R}^N.$ Other asymptotics of the minimal speed were proved in El Smaily \cite{El Smaily} and El Smaily, Hamel, Roques \cite{EHR}.

In the case where $N=2,$ we give necessary and sufficient conditions on the streamlines (or the trajectories) of the advection field $q$ for which
the limit (\ref{large advection}) is positive or null.
In the following lemma, we describe the family of ``unbounded periodic trajectories'' of $L-$periodic 2-dimensional vector fields $q$.
\vskip0.25cm
\begin{lemma}[unbounded periodic trajectories]\label{periodic trajs}
We assume here that $N=2$ and hence $d\in\{1,2\}$. Assume that $q$ satisfies (\ref{cq}). Let $T(x)$ be an unbounded periodic trajectory of $q$ in $\Omega.$ That is, there exists $\mathbf{a} \in L_1\Z\times L_2\Z\setminus\{0\}$ (resp. $L_1\Z\times \{0\}\setminus\{0\}$) when $d=2$ (resp. $d=1$) such that $T(x)=T(x)+\mathbf{a}.$ In this case, we say that $T(x)$ is $\mathbf{a}-$periodic. Then, if $T(y)$ is another unbounded periodic trajectory of $q,$ $T(y)$ is also $\mathbf{a}-$periodic.

Moreover, in the case $d=1,$ $\mathbf{a}=L_1e_1.$ That is, all the unbounded periodic trajectories of $q$ in $\Omega$ are $L_1e_1-$periodic.
\end{lemma}
\vskip0.2cm
\begin{theorem}\label{case N=2}
Assume that $N=2$ and that $\Omega$ and $q$ satisfy (\ref{comega}) and (\ref{cq}) respectively. The two following statements are equivalent:

(i) There exists $w \in \mathcal{I},$ such that $\displaystyle \int_C qw^2 \neq 0$.

(ii) There exists a periodic unbounded trajectory $T(x)$ of $q$ in $\Omega$.

Moreover, if (ii) is verified and $T(x)$ is $\mathbf{a}-$periodic, then for any $w\in\mathcal{I}$ we have $\ds{\int_Cq\,w^2\in\R\mathbf{a}}.$
\end{theorem}

\begin{remark}
The periodicity assumption on the trajectory in (ii) is crucial. Indeed there may exist \\unbounded trajectories which are not periodic, even though the vector field $q$ is periodic. Consider, for example,the following function $\phi$:
$$
\phi(x,y) :=\left\{\begin{array}{l}\ds{ e^{-\frac{1}{\sin^2(\pi y)}}\sin(2\pi(x+\ln(y-[y])))\text{ if }y\not\in\mathbb{Z}},\vspace{3pt}\\
0\text{ otherwise},
\end{array}\right.
$$
where $[y]$ denotes the integer part of $y$. This function is $C^\infty$ on $\mathbb{R}^2$, and $1$-periodic in $x$ and $y$. Hence the vector field
$$
q = \nabla^\perp \phi
$$
is also $C^\infty$, $1$-periodic in $x$ and $y$, and verifies $\int_{[0,1]\times[0,1]} q = 0$ with $\nabla\cdot q\equiv0.$ A quick study of this vector field shows that the part of the graph of $x \mapsto e^{-x}$ lying between $y=0$ and $y=1$ is a trajectory of $q$, and is obviously unbounded and not periodic. Moreover, there exist no periodic unbounded trajectories for this vector field, so the theorem asserts that  for all $w \in \mathcal{I}$ we have
$$
\int_C qw^2 =0.
$$
\end{remark}

\begin{corollary}\label{NSC for lim to be zero} Assume that $N=2$ and that $\Omega,$ $A,$ $q$ and $f$ satisfy the conditions (\ref{comega}), (\ref{cA}), (\ref{cq}) and (\ref{cf1}-\ref{cf2}) respectively. Then,

(i) If there exists no periodic unbounded trajectory of $q$ in $\Omega,$ then
$$
\lim_{M\rightarrow+\infty}\ds{\frac{\ds{c^{*}_{\Omega,A,M\,q,f}(e)}}{M}}=0,
$$ for any unit direction $e.$

(ii) If there exists a periodic unbounded trajectory $T(x)$ of $q$ in $\Omega$ {\rm(}which will be $\mathbf{a}-$periodic for some vector $\mathbf{a}\in\mathbb{R}^2${\rm)} then
\begin{equation}\label{null equivalence}
\lim_{M\rightarrow+\infty}\ds{\frac{\ds{c^{*}_{\Omega,A,M\,q,f}(e)}}{M}} >0\,\Longleftrightarrow\,\tilde{e}\cdot\mathbf{a}\neq0.
\end{equation}
We mention that in the case where $d=1,$ we have $\tilde{e}=\pm e_1.$  Lemma \ref{periodic trajs} yields that $\tilde{e}\cdot\mathbf{a}=\pm L_1\neq0.$ Referring to (\ref{null equivalence}), we can then write,  for $d=1,$
\begin{equation}\label{null equivalence d=1}
\lim_{M\rightarrow+\infty}\ds{\frac{\ds{c^{*}_{\Omega,A,M\,q,f}(e)}}{M}} >0\,\Longleftrightarrow\,\text{(there exists a periodic unbounded trajectory $T(x)$ of $q$ in $\Omega$)}.
\end{equation}

\end{corollary}

In order to prove Theorem \ref{case N=2}, we proved the following qualitative property for the family of the advection fields that we consider in
this work. The following proposition can be viewed as a generalization of the Hodge representation to unbounded domains with a periodic structure.
\vskip0.2cm
\begin{proposition}\label{phodge}
Let $N=2,$ $d=1$ or $2$ where $d$ is defined in (\ref{comega}). Let $q \in C^{1,\delta}(\overline{\Omega})$, $L$-periodic with respect to $x$ and verifying the conditions (\ref{cq}).
Then, there exists $\phi \in C^{2,\delta}(\overline{\Omega})$, $L$-periodic with respect to $x$, such that
\begin{equation}\label{hodge}
q = \nabla^{\perp} \phi \ \text{ in } \Omega.
\end{equation}
Moreover, $\phi$ is constant on every connected component of $\partial \Omega$.
\end{proposition}
\vskip0.2cm
\begin{remark} We mention that the representation $q = \nabla^{\perp} \phi$  is already known in the case where the domain $\Omega$ is bounded and \textbf{simply connected} or  equal to whole space $\R^2.$ However, the above proposition applies in more cases due to the condition $q\cdot\nu=0$ on $\partial\Omega$. For example, it applies when $\Omega$ is the whole space  $\R^2$ with a periodic array of holes or when $\Omega$ is  an infinite cylinder which may have an oscillating boundary and/or a periodic array of holes.
\end{remark}

\vskip0.3cm

The proofs of all the above results, further details and more asymptotic properties of the KPP minimal speed are shown in \cite{ek1}.


\end{document}